\documentclass{SUG} 
\usepackage{amsfonts,amssymb}

\newcommand{\y}{\\[0pt]}
\newcommand{\q}{\quad}\newcommand{\qq}{\qquad}
\newcommand\vs{\vskip10pt}\newcommand\bb{\bigbreak}
\newcommand\n{\noindent}
\newcommand\ip{\raise1pt\hbox{\large$\lrcorner$}\,}
\newcommand{\hb}[1]{\hbox{\small\q#1\q}}
\newcommand\lra{\longleftrightarrow}
\newcommand\mb{\medbreak}\renewcommand\sb{\smallbreak}
\newcommand{\pd}[1]{\hbox{$\ds\frac{\partial#1}{\partial t}$}}

\newcommand{\ds}{\displaystyle}\newcommand{\ts}{\textstyle}
\newcommand{\g}{\mathfrak{g}}\newcommand{\fh}{\mathfrak{h}}
\newcommand{\ga}{\gamma}\newcommand\w{\omega}
\renewcommand{\vp}{\varphi}\newcommand{\svp}{*\kern-.5pt\vp}
\newcommand{\RP}{\mathbb{RP}}\newcommand{\CP}{\mathbb{CP}}
\newcommand{\C}{\mathbb{C}}\newcommand{\R}{\mathbb{R}}
\newcommand{\B}{\mathbb{B}}\newcommand{\M}{\mathbb{M}}
\newcommand{\T}{\mathbb{T}}
\newcommand{\+}{\!+\!}
\renewcommand{\-}{\!-\!}\newcommand{\e}{\mathrm{e}}
\newcommand{\na}{\nabla}\renewcommand{\a}{\alpha}
\newcommand{\si}{\sigma}\newcommand{\os}{\ol\sigma}
\newcommand{\La}{\Lambda}\newcommand{\ST}{S^2_0\T}
\newcommand{\psm}{{\psi_-}}\newcommand{\psp}{{\psi_+}}
\newcommand{\hd}{\hat d}\newcommand{\hw}{\hat\w}
\newcommand{\hpsp}{\hat\psi_+}\newcommand{\hpsm}{\hat\psi_-}
\newcommand{\cW}{\mathcal{W}}\newcommand{\cX}{\mathcal{X}}
\newcommand{\we}{\wedge}
\newcommand{\op}{\oplus}\newcommand{\ot}{\otimes}
\newcommand{\ol}{\overline}
\newcommand{\rf}[1]{(\ref{#1})}\def\lie#1({\mathfrak{#1}(}
\newcommand{\su}{\mathfrak{su}}
\newcommand{\ba}{\begin{array}}\newcommand{\ea}{\end{array}}
\renewcommand{\ee}[1]{\label{#1}\end{equation}} 
\newcommand\qed{\hfill$\square$\mb}
\def\E{\raise1pt\hbox{$\bigwedge$}}
\newcommand{\ft}[2]{\hbox{$\textstyle\frac{#1}{#2}$}}
\newcommand{\nbf}[1]{\stepcounter{enumii}
                                  \bb\n\textbf{\theenumii\ #1\ }}
\newcommand{\nit}[1]{\sb\n\textit{#1.}}
\renewcommand{\theenumii}{\arabic{section}.\arabic{enumii}}
\begin{document}\parskip1pt\parindent15pt

\title{\LARGE\bf The intrinsic torsion of SU(3) and G\boldmath{$_2$}
structures}

\author{Simon Chiossi}

\address{Dipartimento di Matematica, Universit\`a di Genova, Via Dodecaneso
35,\\ I\,--\,16146 Genova}

\author{Simon Salamon}

\address{Dipartimento di Matematica, Politecnico di Torino, Corso Duca degli
Abruzzi 24,\\ I\,--\,10129 Torino}

\maketitle

\centerline{\it To Antonio Naveira on the occasion of his 60th birthday}\vs

\abstracts{We analyse the relationship between the components of the intrinsic
torsion of an $SU(3)$-structure on a 6-manifold and a $G_2$-structure on a
7-manifold. Various examples illustrate the type of $SU(3)$-structure that can
arise as a reduction of a metric with holonomy $G_2$.}

\section*{Introduction} 

Let $G$ be a subgroup of $SO(n)$. A $G$-structure on a smooth manifold $M$ of
dimension $n$ induces a Riemannian metric $g$ on $M$. The failure of the
holonomy group of the Levi-Civita connection of $g$ to reduce to $G$ is
measured by the so-called intrinsic torsion $\tau$. It is known\cite{Bry,Shol}
that the latter is a tensor which takes values at each point in
$T^*\ot\g^\perp$ where $T^*$ is the cotangent space and $\g^\perp$ is the
orthogonal complement of $\g$ in $\lie so(n)\cong\E^2T^*$.

This note is concerned with the cases \[\ba{l} \hbox{(1)}\q SU(3)\subset
SO(6),\y\hbox{(2)}\q G_2 \subset SO(7).\ea\] The respective tensors $\tau_1$
and $\tau_2$ belong to spaces of dimension 42 and 49. The fact that $SU(3)$ is
a maximal subgroup of $G_2$ gives a direct relationship between the two
structures. Indeed, the sets of reductions (1) and (2) are both parametrized
by the projective space \be\RP^7=\frac{SO(6)}{SU(3)}=\frac{SO(7)}{G_2}.\ee{RP}
The fact that this space itself admits homogeneous $G_2$-structures has
applications to the study of families of $G_2$-structures. Moreover, the
fibration $\RP^7\to \CP^3$ is indicative of the way in which $G_2$-structures
can in general be built from almost-Hermitian structures on a 6-manifold.

We begin by describing the tensor $\tau_1$ determined by an $SU(3)$-structure
on a 6-manifold $M$, thereby refining the theory for $U(3)$. An additional
summand in the $SU(3)$ case can be used to construct a new conformally
invariant torsion tensor. It is well known that a holonomy reduction to $SU(3)$
is characterized by the existence of a symplectic form together with a closed
form of `type $(3,0)$', and it follows that all the components of $\tau_1$ can
be calculated in terms of \textit{exterior} derivatives of the forms defining
the reduction. The special relevance of 3-forms in describing 6-dimensional
structures is already documented,\cite{Hjdg} and this paper presents same
additional applications.

In the general set-up, the $SU(3)$ reduction leads to a splitting of the
Nijenhuis tensor in two equal parts, which give rise to different components
of the tensor $\tau_2$ on a 7-manifold $\M$ whose structure reduces from $G_2$
to $SU(3)$. This aspect of the theory is reminiscent of self-duality in four
dimensions, and the $G_2$ examples analysed in subsequent sections can by
analogy be divided into those of self-dual and anti-self-dual type. The
distinction arises from whether the 7-manifold is foliated by leaves of
dimension 1 or 6.

Our first descriptions of $\tau_1,\tau_2$ are `static' in the sense that they
relate to a fixed $G$-structure and are purely algebraic. We subsequently
examine how the components of $\tau_1$ determine those of $\tau_2$ in various
situations in which the geometry of the 6- and 7-manifolds are interrelated,
with the inclusion $SU(3)\subset G_2$ varying from point to point. The
evolution equations discussed by Hitchin\cite{Hit} are interpreted using the
notion of a half-flat $SU(3)$ structure. We provide additional examples of
incomplete metrics with holonomy $G_2$ of the type discovered by Gibbons et
al\cite{GLPS} that suggest that half-flat structures occur naturally on
6-dimensional nilmanifolds.

The final section undertakes an investigation of certain cases in which the
$G_2$-manifold $\M$ is a circle bundle over a 6-manifold $M$ endowed with an
appropriate structure. We provide an explicit description of $\tau_2$ as a
function of $\tau_1$ and a curvature 2-form, and consider the case of the
canonical circle bundle over a K\"ahler 3-fold. When the holonomy of $\M$
reduces to $G_2$, the quotient $M=\M/S^1$ is a symplectic manifold with a type
of generalized Calabi-Yau geometry that we describe briefly in terms of
$\tau_1$. Examples of such quotients of the known complete metrics with $G_2$
holonomy incorporate interesting global features,\cite{BS,AW} and it is hoped
that the techniques of this paper will aid a fuller investigation of this
situation.

\section{Static SU(3) structures}\setcounter{enumii}0

Let $M$ be a 6-manifold with a $U(3)$-structure. Thus $M$ is equipped with a
Riemannian metric $g$, an orthogonal almost-complex structure $J$ and an
associated 2-form $\w$. The exterior forms on $M$ may be decomposed into types
relative to $J$, and we adopt the following notation\cite{FFS} at each point:
\be\ba{rclcl} T^* &=& [\![\La^{1,0}]\!],&&\y 
\E^2T^* &=&[\![\La^{2,0}]\!]\op[\La^{1,1}]\ 
        &\cong& [\![\La^{2,0}]\!]\op[\La^{1,1}_0]\op\R,\y
\E^3T^* &=& [\![\La^{3,0}]\!]\op[\![\La^{2,1}]\!]
        &\cong&[\![\La^{3,0}]\!]\op[\![\La_0^{2,1}]\!]\op[\![\La^{1,0}]\!],\y
\E^4T^* &=& [\![\La^{3,1}]\!]\op[\La^{2,2}]
        &\cong&[\![\La^{2,0}]\!]\op[\La_0^{1,1}]\op\R.
\ea\ee{ext6}\sb

The induced metric distinguishes the circle $B$ consisting of elements of unit
norm in the 2-dimensional space $[\![\La^{3,0}]\!]$. An $SU(3)$-structure is
determined by the choice of a real 3-form $\psp$ lying in $B$ at each point, or
equivalently a section of the associated $S^1$-bundle $\B$. The associated
$(3,0)$-form is \be\Psi=2(\psp)^{3,0}=\psp+i\psm\ee{Psi} with $\psm=J\psp$. We
may then write \be\B=\{a\psp+b\psm: a^2+b^2=1\},\ee{B} and this description
remains valid locally even if a global reduction from $U(3)$ to $SU(3)$ does
not exist.

To be more explicit, we may choose a local orthonormal basis $(e^1,\ldots,e^6)$
of $T^*$ such that \[ \Psi=(e^1+ie^2)\we(e^3+ie^4)\we(e^5+ie^6).\] Consequently
\be\ba{rcl} \w &=& e^{12}+e^{34}+e^{56},\y \psp &=& e^{135}-e^{146}-e^{236}-
e^{245},\y\psm &=& e^{136}+e^{145}+e^{235}-e^{246},\ea\ee{wpm} where $e^{135}$
stands for $e^1\we e^3\we e^5$ etc. These forms are subject to the
compatibility relations \be\ba{rcl}\w\we\psi_\pm&=&0,\y\psp\we\psm&=&\ft23
\w^3.\ea \ee{compat}

The intrinsic torsion of the $U(3)$-structure can be identified with $\na J$
or $\na\w$ and belongs to the space \[T^*\ot\lie u(3)^\perp=\cW_1\op\cW_2\op
\cW_3 \op\cW_4,\] whose four components were first described by
Gray--Hervella.\cite{GH} The intrinsic torsion $\tau_1$ of the
$SU(3)$-structure lies in the enlarged space \[T^*\ot\su(3)^\perp\cong
T^*\ot([\![\La^{2,0}]\!]\op\R)\cong T\ot(T\op\R),\] given that now
$\La^{2,0}\cong\La^{0,1}$. Thus \[\tau_1\in\cW_1\op\cW_2\op
\cW_3\op\cW_4\op\cW_5,\] with $\cW_5\cong T$. Properties of the various
components are indicated by Table~1.\mb

\begin{table}[t]\caption{$SU(3)$ torsion}
\[\ba{|c|c|c|c|c|}\hline\hb{component} & \hb{dim$_\R$}&\hb{$U(3)$-module} &
\multicolumn{2}{|c|}{\hb{$SU(3)$-module}}\\\hline \cW_1 & 2 &
[\![\La^{3,0}]\!] &\q\R\qq& 
\R\\\hline \cW_2 & 16 & [\![V]\!] & \su(3) &
\su(3)\\\hline \cW_3 & 12 & [\![\La_0^{2,1}]\!] &
\multicolumn{2}{|c|}{[\![S^{2,0}]\!]}\\\hline\cW_4 & 6 & T &
\multicolumn{2}{|c|}{T}\\\hline\cW_5 & 6 & T &
\multicolumn{2}{|c|}{T}\\\hline\ea\]\end{table}

We denote the component of $\tau_1$ in $\cW_i$ in a formal way by $W_i$,
though we shall need to supply more precise definitions shortly. It is well
known that the components of $\na J$ in $\cW_1\op\cW_3\op\cW_4$ are determined
by $d\w$, and those in $\cW_1\op\cW_2$ by the Nijenhuis tensor. The choice of
basis $(\psp,\psm)$ of $[\![\La^{3,0}]\!]$ at each point provides an
isomorphism \[\cW_1\op\cW_2=[\![\La^{2,0}\ot\La^{1,0}]\!]=
[\![\La^{3,0}\ot\La^{1,1}]\!]\cong\R^2\ot\lie u(3).\] The rank of $W_1\+W_2$
(or of $W_2$ on its own) in this tensor product is equal to one of 0,1,2.

\nit{Remark} If $M$ is an almost-Hermitian manifold for which the bundle $\B$
is trivial but not actually trivialized, the basis $(\psp,\psm)$ is defined up
to an overall constant action by $S^1$, and the rank of $W_1+W_2$ is a global
invariant. This situation occurs naturally on 6-dimensional Lie groups of the
type considered below.\mb

On a complex manifold, $d\psp$ belongs to $\La^{3,1}\op\La^{1,3}$ at each
point. It follows that the component of $d\psp$ in $\La^{2,2}$ is determined
by the Nijenhuis tensor, and therefore by $W_1\op W_2$. We may define the two
scalar components $W_1^\pm\in\R$ of $\cW_1$ by \[\ba{l} d\psp\we\w=\psp\we
d\w=W_1^+\,\w^3,\y d\psm\we\w=\psm\we d\w=W_1^-\,\w^3,\ea\] where
$\w^3=\w\we\w\we\w$. Similarly, $W_2=W_2^++W_2^-$ in which
\[\ba{l}(d\psp)^{2,2}=W_1^+\w^2+W_2^+\we\w,\y 
        (d\psm)^{2,2}=W_1^-\w^2+W_2^-\we\w,\ea\]
so that $W_2^\pm\in[\La^{1,1}_0]$ are effective $(1,1)$-forms.

Given that \[d\psp-id\psm=d\ol\Psi\in\La^{1,3}\op\La^{2,2},\] the remaining
components of $d\psp,d\psm$ are related by \be (d\psp)^{3,1}=i(d\psm)^{3,1}.
\ee{lemma} It is now clear that the $W_5$-component of $\tau_1$ arises from
\rf{lemma}. In summary, we have

\nbf{Theorem} The five components of $\tau_1$ are determined by $d\w,\,d\psp,
\,d\psm$, in the following manner: \[\ba{rcl} W_1&\lra& (d\w)^{3,0}\y
W_2&\lra&((d\psp)_0^{1,1},\ (d\psm)_0^{1,1})\y W_3&\lra&(d\w)_0^{2,1}\y
W_4&\lra&\w\we d\w\y W_5&\lra&(d\psi_\pm)^{3,1}\ea\] (refer to
\rf{ext6}).\mb

It is significant that $W_4$ and $W_5$ arise from isotypic summands of the
space $T^*\ot\lie su(3)^\perp$. Before moving on to seven dimensions, it is
convenient to give a more precise definition of these components too in order
that they may be compared directly. We shall do this by means of the
contraction \[\ip:\E^kT^*\ot\E^nT^*\to\E^{n-k}T^*\] that exploits the
underlying Riemannian metric, with the convention that $e^{12}\ip
e^{12345}=e^{345}$ etc.

\nbf{Definition} The components of $\tau_1$ in $\cW_4,\cW_5$ are given by
\[\ba{rcl} W_4&=&\ft12\w\ip d\w,\y W_5&=&\ft12\psp\ip d\psp.\ea\]\mb

\n The coefficient of one half is added with the following examples in mind.
Given $\w$ as in \rf{wpm}, suppose that $d\w=\w\we e^1$ and $d\psp=\psp\we
e^1$. Then \be\ba{rcl} W_4&=&\ft12\w\ip(e^{134}+e^{156})=e^1,\y W_5 &=&
\ft12\psp\ip(e^{1236}+e^{1245})=e^1.\ea\ee{e1} Now suppose that 
\[(d\psp)^{3,1}=\Psi\we\os=(\psp+i\psm)\we\os.\] Then $W_5=2(\os+\si)$, whereas
$\ft12\psp\ip d\psm= 2i(\si-\os)$ by \rf{lemma}. It follows that \be \psp\ip
d\psm=J(\psp\ip d\psp),\ee{swap} a useful re-interpretation of \rf{lemma}.

Each of the components $W_1,W_2,W_3$ is at worst re-scaled under a conformal
change of metric \be g\mapsto\e^{2f}g.\ee{conf} This is a consequence of the
fact that none of the corresponding representations in Table~1 is isomorphic to
the cotangent space $T^*$ containing the 1-form $df$. The reduction to $SU(3)$
permits one to define an additional conformally invariant component:

\nbf{Proposition} The tensor $3W_4+2W_5$ is unchanged by \rf{conf}.\mb

\nit{Proof} The transformation \rf{conf} multiplies 1-forms by $\e^f$. Hence
the exterior derivatives of $\w,\psp$ transform as \[\ba{rcl}d\w &\ \mapsto\ 
&d(\e^{2f}\w) =\e^{2f}d\w + 2\e^{2f}df\we\w,\y d\psp& \ \mapsto\ &
d(\e^{3f}\psp)=\e^{3f}d\psp+3\e^{3f}df\we\psp.\ea\] Retaining $\ip$ exclusively
for the contraction relative to the original metric,\[\ba{rcl} W_4 &\ \mapsto\
& W_4 + \w\ip(df\we\w),\y W_5 &\ \mapsto\ & W_5+\ft32\psp\ip(df\we\psp).\ea\]
The final terms may be evaluated by using \rf{e1} with $df\=e^1$, with the
result that they cancel out in the sum $3W_4+2W_5$.\qed

\section{Static G\boldmath{$_2$} structures}\setcounter{enumii}0

We denote by $\T$ the space $\R^7$, regarded as the standard representation of
the exceptional Lie group $G_2$. The latter acts transitively on the sphere
$S^6$ in $\T$, and the stabilizer of a point of $S^6$ is conjugate to a fixed
subgroup $SU(3)$ of $G_2$. The inclusion $SU(3)\subset G_2$ is therefore
characterized by the orthogonal decomposition \be\T=T\op\R.\ee{6+1} We choose
an orthonormal basis $(e^i)$ of $\T^*$ such that $\a=e^7$ annihilates $T$ at
each point.

A $G_2$-structure on a 7-manifold $\M$ is characterized by a `positive generic
3-form' $\vp$. Adopting a canonical form compatible with \rf{wpm}, we set
\be\ba{rcl} \vp &=& \w\we\a+\psp\y &=&
e^{127}+e^{347}+e^{567}+e^{135}-e^{146}-e^{236} -e^{245}.\ea\ee{vp} The basis
$(e^i)$ is orthonormal for the metric determined by the inclusion $G_2\subset
SO(7)$, and allows us to consider \be\ba{rcl}\svp&=&\psm\we\a +\ft12\w^2\y
&=&e^{1367}+e^{1457}+e^{2357}-e^{2467}+e^{3456}+e^{1256} +e^{1234}.\ea\ee{svp}
The structure of a general $G_2$-manifold will not reduce to $SU(3)$, and
these descriptions are only valid pointwise or locally.

The intrinsic torsion space \[\T^*\ot\g_2^\perp=\cX_1\op\cX_2\op\cX_3\op\cX_4
\] has four components of respective dimensions $1,14,27,7$, first described
by Fernandez--Gray.\cite{FG} Various constructions\cite{Bry,BS,GPP,Joy} of
metrics with holonomy equal to $G_2$ are based on the significant fact that
the holonomy reduction is characterized by the simultaneous closure of $\vp$
and $\svp$.

Given that $\T$ represents the tangent space of $\M$, the $G_2$ counterpart of
\rf{ext6} consists of the decompositions \be\ba{rcl} \T^* &\cong&
\T\y\E^2\T^*&\cong& \g_2\op\T\y \E^3\T^* &\cong& \R\op\T\op\ST.\ea\ee{ext7} It
follows that $\g_2^\perp\cong\T$, and the spaces $\cX_1,\cX_2,\cX_3,\cX_4$ are
isomorphic to $\R$, $\g_2$, $\ST$, $\T$ respectively. The components of
$\na\vp$ can be recovered from those of $d\vp$ and $d\svp$ as follows:
\be\ba{rcl} X_1+X_3+X_4 &\lra& d\vp,\y X_1 &\lra& d\vp\we\vp,\y X_2+X_4 &\lra&
d\svp,\y X_4 &\lra& (*d\svp)\we(\svp).\ea\ee{1234} The $G_2$-structure is
called \textit{calibrated} (respectively \textit{cocalibrated}) if $d\vp=0$
(respectively $d\svp=0$), and examples with all possible combinations of
torsion components are known.\cite{F,Cab,CMS}

\begin{table}[t]\caption{$G_2$ torsion}
\[\ba{|c|c|c|c|c|c|c|}\hline\hb{component}&\hb{dim} & \hb{$G_2$-module}
& \multicolumn{4}{|c|}{\hb{$SU(3)$-module}}\\\hline \cX_1 & 1& \ \R\  &
\multicolumn{1}{|c}{}&\multicolumn{3}{l|}{\ \R\ }\\\hline \cX_2 & 14 &\ \g_2\ &
\multicolumn{3}{|r|}{\ T\ } &\ \su(3)\ \\\hline \cX_3 & 27 &\ST&\ 
[\![S^{2,0}]\!]\ &\ \R\ & \ T\  &\ \su(3)\ \\\hline \cX_4 & 7&\ \T\ & 
\multicolumn{2}{|r|}{\ \R\ } &\multicolumn{2}{l|}{\ T\ }\\\hline\ea\]\end{table}

Whilst $\cX_4$ is isomorphic to \rf{6+1}, there are analogous decompositions
of the other spaces in \rf{ext7}, and therefore of $\cX_2$ and $\cX_3$, under
the action of $SU(3)$. These are described in Table~2. In particular, $\ST$
also contains a 1-dimensional $SU(3)$-module, and it follows that
\be\ba{ccc}&\a&\\ &\w&\\ \psp,&\w\we\a,& \psm\\ \psm\we\a,&\w^2,& \psp\we\a\\
&\w^2\we\a& \\ &\w^3& \ea\ee{list} is a list of the exterior forms on $\M$
fixed by $SU(3)$. With reference to the summands in \rf{ext7}, we may assert
that

\nbf{Lemma} \[\left.\ba{r}\psm\in\T\y 3\psp-4\w\we\a\in\ST\ea\right\}
\subset\E^3\T^*,\qq\left.\ba{r}\psp\we\a\in\T\y 3\psm\we\a-2\w^2\in\ST
\ea\right\}\subset\E^4\T^*.\]

\nit{Proof} Given \rf{compat}, $\ga=3\psp-4\w\we\a$ satisfies $\ga\we(\svp)=0$
and is therefore orthogonal to $\vp$. Since $-\psm$ is the contraction of
$\svp$ with the tangent vector dual to $\a$, it lies in the submodule $\T$ of
$\E^3\T^*$. Since $\ga$ is also orthogonal to $\psm$, it lies in $\ST$. The
invariant 4-forms are obtained by observing that $*\psm=\psp\we\a$ and
$*(\w\we\a)=\frac12\w^2$.\qed\sb

\nit{Remark} The existence of various $SU(3)$-invariant elements of $\E^3\T$
gives rise to a choice of induced $G_2$-structures in the passage from 6 to 7
dimensions. For example,
\[3\psp-4\w\we\a=3\left[ \w\we(-\ft43\a)+\psp\right]\] determines a $G_2$
structure with reversed orientation on $\T$ and different scalings relative to
\rf{6+1}. Other choices of coefficients will have the effect of modifying
combinations in Theorem~3.1 below.\sb

The three components of $\T^*\ot\g_2^\perp$ isomorphic to $T$ can be detected
from corresponding components of $d\vp$ and $d\svp$. It is useful to record
the following list for diagnostic purposes.

\nbf{Lemma} \[\ba{rcccll} \zeta &=& e^{1347}+e^{1567}-e^{1236}-e^{1245}&\in
T&\subset\T&\subset\E^4\T^*,\y \eta&=&e^{1347}+e^{1567}+e^{1236}+e^{1245}&\in
T&\subset\ST&\subset\E^4\T^*,\y \xi &=& e^{13456}+e^{12357}-e^{12467}&\in
T&\subset\T&\subset\E^5\T^*,\y \vartheta &=&2e^{13456}-e^{12357}+e^{12467}&\in
T&\subset\g_2&\subset\E^5\T^*.\ea\]

\nit{Proof} Each of these forms represents the element of $T$ dual to $e^1$ in
an appropriate guise. For example, $\zeta=e^1\we\vp$ and $\eta$ is a linear
combination of $\zeta$ and $\psp\we e^1=e^{1236}+e^{1245}$ orthogonal to
$\zeta$. We may define $\xi$ as $e^1\we(\svp)=*(e^1\ip\vp)$. Then $\vartheta$
is a linear combination of $\xi$ and $e^1\we \w^2=2e^{13456}$ orthogonal to
$\xi$.\qed

\section{Product manifolds}\setcounter{enumii}0

We now suppose that $\M$ is a 7-manifold with an $SU(3)$-structure, so that
the differential forms $\a,\w,\psp,\psm$ of respective degrees $1,2,3,3$ and
constant norm are all defined globally. In this and the following sections, we
shall investigate properties of the $G_2$-structure defined with the
convention of \rf{vp}. In general, one may write $d\a=\a\we\beta+\ga$, where
$\beta,\ga$ are forms with values in the subspace $T^*$ at each point. For
example, the equation $\ga=0$ is the integrability condition for the
6-dimensional distribution now determined by \rf{6+1}. We shall consider
various special cases, the simplest of which is that in which $\M$ is the
Riemannian product of $M$ with an interval or circle, so that $\na\a$ (and so
$d\a$) vanishes.

In the product situation, we choose to write $\a=e^7=dt$, so that \be\ba{rcl}
d\vp &=& d\w\we dt+d\psp,\y d\svp &=& d\psm\we dt+\w\we d\w.\ea\ee{dd} Let
\[\ba{rcl} \tau_1 &=& (W_1^+\+W_1^-)+(W_2^+\+W_2^-)+W_3+W_4+W_5,\y\tau_2 &=&
X_1+X_2+X_3+X_4\ea\] denote the respective intrinsic torsion tensors, as
defined in the previous sections. Since $\na\vp$ can be computed in terms of
$\na\w$ and $\na\psp$, the tensor $\tau_2$ is determined by $\tau_1$. Any
$SU(3)$-invariant component of $\tau_2$ must be a linear combination of $W_1^+$
and $W_1^-$, and any component isomorphic to $T$ a linear combination of $W_4$
and $W_5$. The precise statement is

\nbf{Theorem} The four components of $\tau_2$ are determined by the seven 
components of $\tau_1$ as follows. \[\ba{rcl} 
X_1&\lra& W_1^+\y
X_2&\lra& (W_2^-,\ 2W_4\+W_5)\y
X_3&\lra& (W_1^+,\ W_2^+,\ W_3,\ W_4\+W_5)\y
X_4&\lra& (W_4\-W_5,\ W_1^-).\ea\]\sb

\nit{Proof} The component $X_1$ arises from \[\vp\we d\vp=\psp\we d\w\we dt+
\w\we dt\we d\psp=2W_1^+\w^3\we\a.\] Similarly, $X_4$ is determined by the
contraction \[\vp\ip d\vp=-\w\ip d\w+\psp\ip d\psp+(\psp\ip d\w)dt.\] The first
two terms of the right-hand side are $-2W_4$ and $2W_5$ by Definition~1.2, and
the last term corresponds to $d\w\we\psm$ (given that $\psp\we\psp=0$) and so
$W_1^-$. This justifies the description of $X_4$.

The hypotheses $d\w=\w\we e^1$ and $d\psm=\psm\we e^1$ are compatible with the
constraint $W_4=W_5$. In this case \[d\svp=-e^{12357}+e^{12467}+2e^{13456}=
\vartheta,\] confirming that $X_4=0$. In order to obtain $\xi$ instead of
$\vartheta$, we need to take $d\w=\ft12\w\we e^1$ and $d\psm=-\psm\we e^1$,
which corresponds to $2W_4+W_5=0$.

The association of $W_2^-$ with $X_2$ and $W_2^+,W_3$ with $X_3$ follows
immediately from \rf{1234}. The hypotheses $d\w=\w\we e^1$ and $d\psp=-\psp\we
e^1$ are compatible with the constraint $W_4+W_5=0$. This implies that
\[d\vp=e^{1347}+e^{1567}-e^{1236}-e^{1245}=\zeta,\]
whence the $T$-component of $X_3$ is proportional to $W_4+W_5$.\qed

\nit{Remark} The difference $49-42=7$ of the dimensions of the spaces
containing $\tau_1$ and $\tau_2$ is accounted for by the repetition of $W_1^+$
and a linear combination of $W_4,W_5$ in the above list. This redundancy is
eliminated in the more complicated situations described in subsequent
sections. The lack of repetiton between components of $X_1,X_2$ is consistent
with the result\cite{Cab} that a connected $G_2$ manifold with
$\tau_2\in\cX_1\op \cX_2$ has at least one of $X_1,X_2$ zero.

\nbf{Corollary} Suppose that $M$ has an $SU(3)$-structure. The $G_2$-structure
defined on $M\times\R$ by \rf{vp} is cocalibrated if and only if $\tau_1\in
\cW_2^-$.\mb

\nit{Examples. 1} An almost-Hermitian 6-manifold is called
nearly-K\"ahler\cite{Gr} if $\na J$ belongs to the space $\cW_1$. Assuming that
$\na J\ne0$, the structure reduces to $SU(3)$ and we may suppose that
$\tau_1\in\cW_1^-$ with $\psp$ is proportional to $d\w$. The product $M\times
S^1$ then has a $G_2$-structure with $\tau_2\in\cX_4$. Alternatively we may
swap the roles of $\psp,\psm$ to obtain $\tau_2\in\cX_1\op\cX_3$.

\nit{2} A known example\cite{F} of a calibrated nilmanifold can be interpreted
as follows. Let $\g=\R\op\fh$ be a 6-dimensional Lie algebra with structure
determined by \[ de^i=\left\{\ba{ll} 0,&i=1,2,4,5\y e^{25},\qq &i=3,\y
-e^{24},&i=6.\ea\right.\] The definitions \rf{wpm} furnish an associated
nilmanifold $M=\Gamma\backslash G$ with an $SU(3)$-structure for which \[
d\w=0,\q d\psp=0,\q d\psm=e^{1234}-e^{1256},\] whence $\tau_1\in\cW_2^-$. It
follows that $M\times S^1$ has both a calibrated $G_2$ structure and (swapping
$\psp,\psm$) a cocalibrated one with $\tau_2\in\cX_3$. The same Lie algebra
$\g$ was incidentally used\cite{Saps} in the construction of a closed
non-parallel 4-form with stabilizer $Sp(2)Sp(1)$.

\section{Dynamic G\boldmath{$_2$} structures}\setcounter{enumii}0

Let $M$ be a fixed 6-manifold. Suppose that $(\w,\psp,\psm)$ is an $SU(3)$
structure that depends on a real parameter $t$ lying in some interval
$I\subseteq\R$, so that one may regard $\M=M\times I$ as a warped product
fibring over $I$.

To avoid confusion, we denote exterior differentation on $M$ by $\hd$ in this
section. Adopting a unit 1-form $dt$ on $I$ allows us to write \[\ba{rcl} d\vp
&=& (\hd\w-\pd\psp)\we dt+\hd\psp,\y d\svp &=& (\hd\psm+\w\we\pd\w)\we
dt+\w\we\hd\w.\ea\] This motivates

\nbf{Definition} An almost Hermitian 6-manifold is \textit{half-flat} if it
possesses a a reduction to $SU(3)$ for which $d\psp=0$ and $\w\we d\w=0$.\mb

\n Half-flatness is therefore characterized by the closure of $\psp$ and
$\w^2$. It amounts to requiring that $W_1+W_2$ has rank one and that both
$W_4,W_5$ vanish. This eliminates $1+8+6+6=21$ of the total 42 dimensions of
$\tau_1$, which is constrained to lie in $\cW_1^-\op\cW_2^-\op\cW_3$.\sb

If we now suppose that the $G_2$-structure on $\M$ has holonomy group contained
in $G_2$, we may conclude that $M$ is half-flat for all $t$, and that the forms
evolve according to the equations \be\left\{\ba{l} \hd\w=\pd\psp,\y
\hd\psm=-\w\we\pd\w.\ea\right.\ee{evol} Conversely, suppose we are given a
half-flat $SU(3)$-structure at time $t=t_0$. The equations \rf{evol} may be
regarded as a system constraining a closed 3-form $\psp$ and a closed 4-form
$\w^2$. For, if (as in our situation) the stabilizer of $\psp$ is $SL(3,\C)$
then $\psp$ determines $\psm$ via \rf{Psi}. In this way, Hitchin\cite{Hit}
proved that the compatibility equations \rf{compat} are conserved in time, and
this leads to

\nbf{Theorem} Let $M$ be an almost Hermitian 6-manifold which is
half-flat. Then there exists a metric with holonomy contained in $G_2$ on
$M\times I$ for some interval $I$.\mb

A key example of this construction is the following. Given $(M,\hat g)$,
consider the conical metric $g=t^2\hat g+dt^2$ on $M\times\R^+$. Consistent
with this, we set \[\w=t^2\hw,\q \psp=t^3\hpsp,\q\psm=t^3\hpsm,\] where
circumflex indicates a form independent of $t$. Then \rf{evol} becomes
\[d\hw=3\hpsp,\q d\hpsm=-2\hw^2.\] These equations determine the
nearly-K\"ahler class for which $\tau_1\in\cW_1^-$.\sb

\nit{Further examples. 1} We now construct metrics with holonomy $G_2$
associated to each of the three nilmanifolds with $b_1=4$ whose $U(3)$
structures were classified by Abbena et al.\cite{AGS2} Starting with the
Iwasawa manifold $M$, we modify the usual basis of 1-forms in order that \be
de^i= \left\{\ba{ll} 0,&i=1,2,3,4,\y -e^{14}-e^{23},\qq &i=5,\y-e^{13}-e^{42},
&i=6.\ea\right.\ee{iwa} The metric \[\ts g=t^2\sum\limits _{i=1}^4 e^i\ot
e^i+t^{-2}\sum\limits_{i=5}^6 e^i\ot e^i+t^4dt^2\] is compatible with the
natural fibration $M\to T^4$, and the forms
\[\ba{rcl}\w &=& t^2(e^{12}+e^{34})+t^{-2}e^{56},\y \psp+i\psm &=&
t(e^1+ie^2)\we(e^3+ie^4)\we(e^5+ie^6)\ea\] determine a reduction to $SU(3)$ for
which \[\hd\w=t^{-3}\psi_+,\q \hd\psm=-4t\,e^{1234},\] whence $\tau_1\in
\cW_1^-\op\cW_2^-$. Indeed, $g$ is a metric with holonomy $G_2$ on
$M\times\R^+$.

\nit{2} The previous example was discovered by Gibbons et al\cite{GLPS} who
exhibit it as arising from the standard complete metric with holonomy $G_2$ by
`contracting' the isometry group $SO(5)$. A more complicated 2-step example is
associated to the Lie algebra with \[de^i=\left\{\ba{ll} 0,&i=1,2,3,4,\y
-e^{14}-e^{23},\qq & i=5,\y e^{24},&i=6.  \ea\right.\] Consider an orthonormal
basis of 1-forms \[te^1,\q t^2e^2,\q te^3,\q t^2e^4,\q t^{-2}e^5,\q
t^{-1}e^6,\q 2t^4dt,\] and the $SU(3)$-structure for which \[\ba{rcl} \w &=&
t^3(e^{12}+e^{34})+t^{-3}e^{56},\y \psp+i\psm &=&
(e^1+ite^2)\we(e^3+ite^4)\we(e^5+ite^6),\ea\] by analogy to \rf{wpm}. This
yields a $G_2$-structure with closed forms \[\ba{rcl}\vp
&=&2t^7(e^{12}+e^{34})\we dt +2te^{56}\we
dt+e^{135}-t^2(e^{146}+e^{236}+e^{245}),\y \svp &=& -2t^7e^{246}\we
dt+2t^5(e^{145}+e^{136}+ e^{235})\we dt
+e^{1256}+e^{3456}+t^6e^{1234}.\ea\]\mb

\nit{3} The nilpotent 3-step Lie algebra for which \[ de^i=\left\{\ba{ll}
0,&i=1,2,4,5\y e^{25},\q &i=3,\y e^{14}-e^{23},\qq &i=6\ea\right.\] gives rise
to an example satisfying the hypotheses of Theorem~4.2. In fact, the forms
\rf{wpm} satisfy \[\w\we d\w=0,\q d\psp=0,\q d\psm=-e^{1256}.\] An explicit
determination of the metric requires an analysis of spaces of invariant exact
forms.\cite{Hit}

\section{Fibred G\boldmath{$_2$} manifolds}\setcounter{enumii}0

Let $(M,\hat g)$ be a Riemannian 6-manifold. In this final section, we consider
the 7-dimensional total space of a circle fibration $\pi:\M\to M$, endowed with
a metric of the form \be g=\a\ot\a+\pi^*\hat g,\ee{g} with $d\a=\pi^*\rho$ for
some 2-form $\rho$ on $M$.

To begin with, suppose that $M$ has an $SU(3)$-structure. The $G_2$-structure
on $\M$ determined by \rf{vp} satisfies the following enhanced version of
\rf{dd}, in which we omit the pullback operator $\pi^*$: \be\ba{rcl}d\vp &=&
d\w\we\a+d\psp+ \w\we\rho.\y d\svp &=& d\psm\we\a+\w\we
d\w-\psm\we\rho.\ea\ee{ddd} Indicating the components of the 2-form $\rho$ by
\[\ba{ccccccc} \E^2T^* &=& \langle\w\rangle &\op& [\La^{1,1}_0]
&\op&[\![\La^{2,0}]\!]\y \rho &=& \rho_0\w &+& \rho_1 &+& \rho_2,\ea\] we
obtain the following generalization of Theorem 3.1.

\nbf{Theorem} The four components of $\tau_2$ are given by\[\ba{rcl} 
X_1&\lra& W_1^++\rho_0\y
X_2&\lra& (W_2^-,\ 2W_4\+W_5\-2\rho_2)\y
X_3&\lra& (3W_1^+-4\rho_0,\ W_2^+\+\rho_1,\ W_3,\ W_4\+W_5\+\rho_2)\y
X_4&\lra& (W_4\-W_5\-\rho_2,\ W_1^-).\ea\]\sb

\nit{Proof} Whilst $X_1$ corresponds to $\vp\we d\vp$, Lemma~2.1 tells
us that the 1-dimensional component in $X_3$ arises from
$(3\w\we\a-4\psp)\we d\vp$. The various coefficients of $\rho_2$ can
be deduced from the observations: (i) $d\vp=0$ implies $d\w=0$, and
(ii) $d\svp=0$ implies $d\psm=0$.\qed

A simple case is that in which the entire torsion $\tau_1$ vanishes, which
corresponds to \[\ba{rcl} d\vp&=&\w\we\rho,\y d\svp&=&-\psm\we\rho=-\psm\we
\rho_2.\ea\] This situation applies when $M$ is the torus $T^6=
\R^6/\mathbb{Z}^6$ endowed with a constant $SU(3)$-structure. The set \rf{RP}
of such structures compatible with $\hat g$ is isomorphic to the set of
$G_2$-structures compatible with \rf{g} (assuming that orientations are also
preserved). We deduce that $\M$ is calibrated if and only if $\rho=0$, and
cocalibrated if and only if $\rho_2=0$.\mb

\nit{Examples. 1} Let $\M$ be a circle bundle with curvature 2-form \[
\rho=e^{12}+e^{34}+e^{56}\] over $T^6$. Then there are no calibrated $G_2$
structures on $\M$ compatible with $g$. The set $\mathcal{C}^*$ of cocalibrated
structures corresponds to $SU(3)$-structures relative to which $\rho$ has type
$(1,1)$. By considering first the space $\CP^3$ (regarded as the set
$SO(6)/U(3)$ of orthogonal almost complex structures\cite{AGS2,AGI}),
$\mathcal{C}^*$ can be seen to be a disjoint union $\RP^1\sqcup\RP^5$.

\nit{2} We can apply the above theory to the real projective space $\RP^7$ by
representing it as $SO(5)/SU(2)$, which fibres over $SO(5)/U(2)\cong\CP^3$. It
is well known that the latter has a homogeneous nearly-K\"ahler metric, so
that there is a reduction to $SU(3)$ with $\tau_1\in\cW_1^\pm$. Moreover, the
curvature 2-form $\rho$ of the above Hopf fibration is an element of
$[\La^{1,1}_0]$ relative to the non-integrable complex structure.\cite{RC}
Example~1 of Section~3 now implies that $\RP^7$ has a $G_2$-structure with
$\tau_2\in\cX_1\op\cX_3$ and another with $\tau_2\in\cX_3\op\cX_4$.

\nit{3} Let us return to the descriptions \rf{RP}. The space of
$G_2$-invariant differential forms on $\RP^7$ is 1-dimensional, and it follows
that $\RP^7$ has a nearly-parallel $G_2$-structure, one with $\tau_2\in
\cX_1$. More generally, we may regard \rf{list} as a list of $SO(6)$-invariant
differential forms on $\RP^7$ subject to the relations $d\a=-\ft12\w$ and \be
d\psi_\pm=\a\we\psi_\mp,\ee{pm} that also give rise to $G_2$-structures with
$\tau_2\in\cX_1\op\cX_3$.\mb

One can arrive at \rf{pm} by considering the canonical $S^1$-bundle $\B$ over
a K\"ahler manifold $M$ of real dimension 6 (see \rf{B}). A local orthonormal
basis $\{\psi^1,\,\psi^2\}$ of sections of the bundle with fibre
$[\![\La^{3,0}]\!]\cong\R^2$ gives rise to coordinates $a^1,a^2$ and radial
parameter $r=\sqrt{(a^1)^2+(a^2)^2}$ on $\B$. The K\"ahler condition implies
that $d\psi^1=\si\we\psi^2$ for some 1-form $\si$ on $M$, and we use \[ b^1=
da^1-a^2\pi^*\si,\q b^2=da^2+a^1\pi^*\si\] to define global forms \[\ba{c}
rdr=a^1b^1+a^2b^2,\qq\a=a^1b^2-a^2b^1,\y \psp=a^1\psi^1+a^2\psi^2,\qq\psm=
a^1\psi^2-a^2\psi^1.\ea\] A straightforward calculation\cite{Shol} shows that
\[rdr\we\psp+\a\we\psm=r^2(b^1\we\psi^1+b^2\we\psi^2)=r^2d\psi_+,\] so that
restricting to $\B$ ($r=1$) we obtain \rf{pm}. Similarly, $d\a=\pi^*\rho$,
where $\rho=d\si$ can be identified with the Ricci form. This leads to the
following result of Baum et al,\cite{BFGK} which forms part of a more general
theory of nearly-parallel $G_2$ structures.\cite{FKMS,CMS}

\nbf{Theorem} If $M$ is K\"ahler-Einstein with positive scalar curvature then
$\B$ has a $G_2$-structure with $\tau_2\in\cX_1$.\bb

We conclude this article by returning to the initial set-up of this section,
together with the assumption that the holonomy of the metric \rf{g} reduces to
the subgroup $G_2$ determined by the 3-form \rf{vp}. This allows us to set
$X_i=0$ in Theorem~5.1.

\nbf{Corollary} If the holonomy group of the metric \rf{g} reduces to $G_2$
then the quotient $M=\M/S^1$ has an $SU(3)$-structure for which $\tau_1\in
\cW_2^+$.\mb

\n Observe that the resulting condition on $\tau_1$ involves a change of sign
from that in Corollary~3.2. Indeed, from \rf{ddd}, we obtain $d\w=0$ and
$d\psm=0$. It is convenient to regard $SU(3)$-structures with
$(d\psm)^{2,2}=0$ as `self-dual', and those with $(d\psp)^{2,2}=0$ as
`anti-self-dual'. The latter type occurred naturally in Sections 3 and 4, and
we focus attention on the $(2,2)$ components so that the terminology is
conformally invariant. In the present situation, we may therefore say that $M$
has a self-dual symplectic structure with \[d\psp=\w\we\rho,\] where
$\rho=\rho_1$ is a closed effective $(1,1)$-form.

Whilst the incomplete examples of Section~4 admit quotients of this type, a
more realistic generalization of the condition $\tau_1=0$ is obtained by
dropping the assumption that the $S^1$ orbits on $\M$ have constant size. In
this situation, there exists a function $f$ on $M$ for which $d(\e^{2f}\w)
=0$. It follows that, applying the conformal transformation \rf{conf}, $M$
once again has a self-dual symplectic structure, though this time $W_5=-df$ is
non-zero and \[\rho=\rho_1+\rho_2=-e^{-3f}(W_2^++\,2df\ip\psm)\] is closed. A
study of this particular class of structures may be valuable in the
construction of metrics with holonomy $G_2$.\mb

\section*{Acknowledgments} This paper was conceived for the conference
`Differential Geometry Valencia 2001', though the material was developed in
the light of feedback from the Durham-LMS Symposium `Special Structures in
Differential Geometry' that took place a month later. The authors wish to
thank the organizers of both conferences.

\def\newblock{}
\bibliography{SUG}

\enddocument